\documentclass[11pt]{article}
\usepackage{amsmath}
\usepackage{amssymb}
\usepackage{latexsym}
\usepackage{amsthm}
\usepackage{tabularx}
\usepackage{booktabs}
\usepackage{mathrsfs}
\usepackage{enumitem}
\usepackage{ytableau}
\usepackage{graphicx}
\usepackage{algorithm,algorithmic}
\usepackage[colorlinks,
            linkcolor=red,
            anchorcolor=blue,
            citecolor=green]{hyperref}

\usepackage{epic}

\newtheorem{thm}{Theorem}[section]

\newtheorem{cor}[thm]{Corollary}
\newtheorem{conj}[thm]{Conjecture}

\allowdisplaybreaks

\begin{document}
\begin{center}
{\large \bf  Schur positivity and log-concavity related to longest increasing subsequences}
\end{center}

\begin{center}
Alice L.L. Gao$^1$, Matthew H.Y. Xie$^2$ and Arthur L.B. Yang$^{3}$\\[6pt]

Center for Combinatorics, LPMC\\
Nankai University, Tianjin 300071, P. R. China\\[6pt]

Email: $^{1}${\tt gaolulublue@mail.nankai.edu.cn},
       $^{2}${\tt xiehongye@163.com},
       $^{3}${\tt yang@nankai.edu.cn}
\end{center}

\noindent\textbf{Abstract.}
Chen proposed a conjecture on the log-concavity of the generating function for the symmetric group with respect to the length of longest increasing subsequences of permutations. Motivated by Chen's log-concavity conjecture,
B\'{o}na, Lackner and Sagan further studied similar problems by restricting the whole symmetric group to certain of its subsets.
They obtained the log-concavity of the corresponding generating functions for these subsets by using the hook-length formula.
In this paper, we generalize and prove their results by establishing the Schur positivity of certain symmetric functions. This also enables us to propose a new approach to Chen's original conjecture.

\noindent \emph{AMS Classification 2010:} 05A05, 05A20.


\section{Introduction}\label{sec-schur_intro}

Given {positive integers} $m,n$ and $\lceil\frac{n}{2}\rceil\leq k \leq n$, let $(k^m,(n-k)^m)$ denote the partition with $m$ parts equal to $k$
and $m$ parts equal to $n-k$. Similarly, for $1\leq k \leq n$, let $(k^m,1^{m(n-k)})$ denote the partition with $m$ parts equal to $k$
and $m(n-k)$ parts equal to $1$. Given a partition $\lambda$, let $f^{\lambda}$ denote the number of standard Young tableaux of shape $\lambda$. The main objective of this paper is to
prove the following result.

\begin{thm}\label{thm-bona-lackner-sagan}
Suppose that $m,n$ are two positive integers.
\begin{itemize}
\item[(1)]
For $\lceil\frac{n}{2}\rceil< k < n$ we have
\begin{align*}
(f^{(k^m,(n-k)^m)})^2&\geq f^{((k+1)^m,(n-k-1)^m)}f^{((k-1)^m,(n-k+1)^m)}.
\end{align*}

\item[(2)]
For $1< k < n$ we have
\begin{align*}
(f^{(k^m,1^{m(n-k)})})^2&\geq f^{((k+1)^m,1^{m(n-k-1)})}f^{((k-1)^m,1^{m(n-k+1)})}.
\end{align*}
\end{itemize}
\end{thm}

The roots of this paper lie in the work by B\'{o}na, Lackner and Sagan \cite{bona2016aoc}, who first proved the above theorem for the case of $m=1,2$ by using the celebrated hook-length formula. We will present two proofs of Theorem \ref{thm-bona-lackner-sagan}, one of which is the same as B\'{o}na, Lackner and Sagan's proof for small $m$, and the other is based on some results on Schur positivity due to Lam, Postnikov, and Pylyavskyy \cite{lam2007ajm}.

Let us first review some backgrounds. We will adopt the notation and terminology found in B\'{o}na, Lackner and Sagan \cite{bona2016aoc}.
Given a positive integer $n$, let $\mathfrak{S}_n$ be the symmetric group of all permutations of $[n]:=\{1,2,\ldots,n\}$. For a given permutation $\pi\in \mathfrak{S}_n$, let $\ell(\pi)$ denote the length of a longest increasing subsequence of $\pi$. Define $L_{n,k}$ to be the set of permutations $\pi\in\mathfrak{S}_n$ with $\ell(\pi)=k$ for $1\leq k \leq n$. Let $\ell_{n,k}=|L_{n,k}|$.  Chen proposed the following conjecture.

\begin{conj}[{\cite[Conjecture 1.1]{chen2008}}]\label{chen_logcave_permu}
For any fixed $n$,
the sequence $\{\ell_{n,k}\}_{k=1}^n$ is log-concave, namely,
$\ell_{n,k}^2\geq \ell_{n,k+1}\ell_{n,k-1}$ for $1<k<n$.
\end{conj}

B\'{o}na, Lackner and Sagan \cite{bona2016aoc} further made a companion conjecture for involutions. Define $I_{n,k}$ to be the set of involutions $\pi\in\mathfrak{S}_n$ with $\ell(\pi)=k$ for $1\leq k \leq n$. Let  $i_{n,k}=|I_{n,k}|$. They proposed the following conjecture.

\begin{conj}[{\cite[Conjecture 1.2]{bona2016aoc}}]\label{bona_logcave_permu}
For any fixed $n$,
the sequence $\{i_{n,k}\}_{k=1}^n$ is log-concave.
\end{conj}

B\'{o}na, Lackner and Sagan showed that there is a close connection between Conjecture \ref{chen_logcave_permu} and Conjecture \ref{bona_logcave_permu} by using the Robinson-Schensted correspondence.
It is well known that each permutation $\pi\in\mathfrak{S}_n$, under the Robinson-Schensted correspondence, is mapped to a pair of standard Young tableaux of the same partition shape, say $\lambda=(\lambda_1,\lambda_2,\ldots)\vdash n$. Moreover, there holds
$\ell(\pi)=\lambda_1$. In that case, we also say that $\pi$ is of shape $\lambda$, denoted $\mathrm{sh}\, \pi=\lambda$. B\'{o}na, Lackner and Sagan proved that if there is a shape-preserving injection from $I_{n,k-1}\times I_{n,k+1}$ to $I_{n,k}\times I_{n,k}$, then there is a shape-preserving injection from $L_{n,k-1}\times L_{n,k+1}$ to $L_{n,k}\times L_{n,k}$, see \cite[Theorem 2.2]{bona2016aoc}.

Though they could not prove Conjectures \ref{chen_logcave_permu} and \ref{bona_logcave_permu}, B\'{o}na, Lackner and Sagan proposed a new way to look at these problems. Given a set $\Lambda$ of partitions of $n$, for $1\leq k \leq n$ let
\begin{align*}
L^{\Lambda}_{n,k}=\{\pi\in L_{n,k}\mid \mathrm{sh}\, \pi\in\Lambda\},\quad \ell^{\Lambda}_{n,k}=|L^{\Lambda}_{n,k}|;\qquad
I^{\Lambda}_{n,k}=\{\pi\in I_{n,k}\mid \mathrm{sh}\, \pi\in\Lambda\}, \quad i^{\Lambda}_{n,k}=|I^{\Lambda}_{n,k}|.
\end{align*}
Thus, the sequence $\{\ell_{n,k}\}_{k=1}^n$ (resp. $\{i_{n,k}\}_{k=1}^n$) is just $\{\ell^{\Lambda}_{n,k}\}_{k=1}^n$ (resp. $\{i^{\Lambda}_{n,k}\}_{k=1}^n$) when taking $\Lambda$ to be the set of all partitions of $n$.
They noted that the log-concavity of $\{\ell^{\Lambda}_{n,k}\}_{k=1}^{n}$ is equivalent to that of $\{i^{\Lambda}_{n,k}\}_{k=1}^{n}$ provided that the set $\Lambda$ contains at most one partition with first row of length $k$ for each $1\leq k \leq n$. They further obtained the following results, see \cite[Theorems 3.1, 3.2, 4.4 and 4.5]{bona2016aoc}.

\begin{thm} \label{thm-bona} Suppose that $n$ is a positive integer and $m=1,2$.
\begin{itemize}
\item[(1)] For $\Lambda=\{(j^m,(n-j)^m) \mid \lceil\frac{n}{2}\rceil\leq j\leq n\}$, the sequence $\{i^{\Lambda}_{mn,k}\}_{k=1}^{mn}$ is log-concave.
\item[(2)] For $\Lambda=\{(j^m,1^{m(n-j)}) \mid 1\leq j\leq n\}$, the sequence $\{i^{\Lambda}_{mn,k}\}_{k=1}^{mn}$ is log-concave.
\end{itemize}
\end{thm}
The Robinson-Schensted correspondence tells that $f^{\lambda}=|\{\pi \mid  \pi^2=id \mbox{ and } \mathrm{sh}\, \pi=\lambda\}|$. Thus Theorem \ref{thm-bona-lackner-sagan} and Theorem \ref{thm-bona} are equivalent to each other for $m=1,2$.

To prove the inequalities on $f^{\lambda}$, a natural way is to use the hook-length formula, as B\'{o}na, Lackner and Sagan did in their paper \cite{bona2016aoc}. Here we will propose another way based on the property of the exponential specialization. Let $\Lambda_{\mathbb{Q}}$ denote the ring of symmetric functions over the field $\mathbb{Q}$ of rational numbers. Recall that the exponential specialization $\mathrm{ex}: \Lambda_{\mathbb{Q}}\longrightarrow \mathbb{Q}[t]$ is defined by acting on the power sums $p_n$ as
\begin{align}\label{exponential_specialization}
\mathrm{ex}(p_n)=t\delta_{1n},
\end{align}
and then extended algebraically. For any symmetric function $f$, let $\mathrm{ex}_1(f)=\mathrm{ex}(f)_{t=1}$. It is well known that
\begin{align}\label{exponential_specialization_s}
\mathrm{ex}_1(s_{\lambda})=\frac{f^{\lambda}}{n!},\quad \mbox{ or equivalently,} \quad f^{\lambda}=\mathrm{ex}_1({n!}s_{\lambda})
\end{align}
for any $\lambda\vdash n$. For more information on the exponential specialization, see \cite{stanlay1999cup2}. Since $\mathrm{ex}$ is an algebra homomorphism, the inequalities on $f^{\lambda}$ considered in Theorem \ref{thm-bona-lackner-sagan} can be deduced from the Schur positivity of the differences of products of Schur functions $s_{\lambda}$.

The rest of the paper is organized as follows.
In Section \ref{sec-schur_2}, we give a proof of
Theorem \ref{thm-bona-lackner-sagan} by using the hook-length formula.
In Section \ref{sec-schur_3}, we present an alternative  proof of
Theorem \ref{thm-bona-lackner-sagan} based on the Schur positivity of certain symmetric functions.

\section{Proof by the hook-length formula}\label{sec-schur_2}

The aim of this section is to give an proof of Theorem \ref{thm-bona-lackner-sagan} by using the hook-length formula.

Let us first give an overview of related definitions and results.
Given a partition $\lambda$, let $\ell(\lambda)$ denote the number of its nonzero parts. Each partition $\lambda$ is associated to a left justified array of cells with $\lambda_i$ cells in the $i$-th row, called the Ferrers or Young diagram of $\lambda$. Here we number the rows from top to bottom and the columns from left to right.
The cell in the $i$-th row and $j$-th column is denoted by $(i, j)$.
The hook-length of $(i,j)$, denoted by $h(i,j)$, is defined
to be the number of cells directly to the right or directly below $(i,j)$, counting $(i,j)$ itself once. The classical hook-length formula is stated as follows, which was discovered by Frame, Robinson and Thrall \cite{frame1954cjm}.

\begin{thm}[\cite{frame1954cjm}]\label{hook_length_schur1}
For any partition $\lambda\vdash n$, we have
\begin{align}\label{hook_length_schur_equa_1}
f^{\lambda}=\frac{n!}{\prod_{(i,j)\in \lambda}h_{(i,j)}}.
\end{align}
\end{thm}
For our purpose here, it turns out to be easier to work with the following equivalent form of the hook-length formula.

\begin{thm}[\cite{fulton1991ny}]\label{hook_length_schur2}
Given a partition $\lambda=(\lambda_{1},\lambda_{2},\cdots,\lambda_{\ell(\lambda)})\vdash n$, we have
\begin{align}\label{hook_length_schur_equa_2}
f^{\lambda}=\frac{n!}{h_{(1,1)}! h_{(2,1)}!\ldots h_{(\ell(\lambda),1)}!}
\prod_{1\leq j_1<j_2\leq \ell(\lambda)}(h_{(j_1,1)}-h_{(j_2,1)}).
\end{align}
\end{thm}

The equivalence between these two formulas is evident by virtue of the equality
\begin{align}\label{hook_length_schur_equa}
\prod_{(i,j)\in \lambda}h_{(i,j)}=
\frac{h_{(1,1)}! h_{(2,1)}!\ldots h_{(\ell(\lambda),1)}!}
{\prod_{1\leq j_1<j_2\leq \ell(\lambda)}(h_{(j_1,1)}-h_{(j_2,1)})}.
\end{align}
It should be mentioned that \eqref{hook_length_schur_equa_2} can be taken as a direct consequence of the Frobenius character formula, see Fulton and  Harris \cite{fulton1991ny}.

Now we can give a proof of Theorem \ref{thm-bona-lackner-sagan}.

\noindent \textit{Proof of Theorem \ref{thm-bona-lackner-sagan}.}
Let us first prove that for $\lceil\frac{n}{2}\rceil< k < n$
\begin{align*}
(f^{(k^m,(n-k)^m)})^2&\geq f^{((k+1)^m,(n-k-1)^m)}f^{((k-1)^m,(n-k+1)^m)}.
\end{align*}
To this end, we will use the expression of $f^{(k^m,(n-k)^m)}$ given by
Theorem \ref{hook_length_schur2}. It is readily to see that the hook-lengths of the first column of the partition $(k^m,(n-k)^m)$ are given by \begin{align}
h_{(i,1)}=\left\{
\begin{array}{ll}
k+2m-i, & \mbox{for $1\leq i\leq m$;}\\
n-k+2m-i, & \mbox{for $m+1\leq i\leq 2m$.}
\end{array}
\right.\label{eq-hooklengths}
\end{align}
Therefore,
\begin{align*}
f^{(k^m,(n-k)^m)}=&\frac{(mn)!}{\prod_{i=1}^m h_{(i,1)}!\times\prod_{i=m+1}^{2m} h_{(i,1)}!}\times\prod_{{{1\leq i_1\leq m}\atop {m+1\leq i_2\leq 2m}}}(h_{(i_1,1)}-h_{(i_2,1)})\\
&\times \prod_{1\leq i_1<i_2\leq m}(h_{(i_1,1)}-h_{(i_2,1)})\times \prod_{m+1\leq i_1<i_2\leq 2m}(h_{(i_1,1)}-h_{(i_2,1)}).
\end{align*}
Substituting \eqref{eq-hooklengths} into the above formula, we obtain
\begin{align*}
f^{(k^m,(n-k)^m)}=&\frac{(mn)!}{\displaystyle{\prod_{i=0}^{m-1}(m+k+i)!\times\prod_{i=0}^{m-1}(n-k+i)!}}
\times
\prod_{0\leq i_1,i_2\leq m-1}(2k-n+1+i_2+i_1)\\
&\times\prod_{0\leq i_1<i_2\leq m-1}(i_2-i_1)\times \prod_{0\leq i_1<i_2\leq m-1}(i_2-i_1).
\end{align*}
Note that the last two factors on the right hand side are independent of $k$. Denote the second factor by $a_k$, namely,
$$a_k=\prod_{0\leq i_1,i_2\leq m-1}(2k-n+1+i_2+i_1).$$
It is easy to verify that
\begin{align*}
\frac{a_k^2}{a_{k-1}a_{k+1}}=&\prod_{0\leq i_1,i_2\leq m-1}
\frac{(2k-n+1+i_2+i_1)^2}{(2k-n-1+i_2+i_1)(2k-n+3+i_2+i_1)}\geq 1,
\end{align*}
since, for any $\lceil\frac{n}{2}\rceil< k <n$ and $0\leq i_1,i_2\leq m-1$, there holds
$${(2k-n+1+i_2+i_1)^2}\geq {(2k-n-1+i_2+i_1)(2k-n+3+i_2+i_1)}$$
 by the inequality of arithmetic and geometric means.
Thus, for $\lceil\frac{n}{2}\rceil< k <n$, we have
\begin{align*}
\frac{(f^{(k^m,(n-k)^m)})^2}{f^{((k+1)^m,(n-k-1)^m)}f^{((k-1)^m,(n-k+1)^m)}}=
\frac{2m+k}{m+k}\times
\frac{n-k+m}{n-k}\times
\frac{a_k^2}{a_{k-1}a_{k+1}}\geq 1,
\end{align*}
as desired.

We proceed to prove the second part of the theorem, namely, for $1< k < n$,
\begin{align*}
(f^{(k^m,1^{m(n-k)})})^2&\geq f^{((k+1)^m,1^{m(n-k-1)})}f^{((k-1)^m,1^{m(n-k+1)})}.
\end{align*}
Let us first give an expression of $f^{(k^m,1^{m(n-k)})}$
by using Theorem \ref{hook_length_schur1}.
Note that, for $1\leq k \leq n$, the hook-lengths of the partition $(k^m,1^{m(n-k)})$ are given by
\begin{align}
     h_{(i,j)}=\left\{\begin{array}{ll}
    m(n-k)+m+k-i, & \mbox{for $j=1$ and $1\leq i\leq m$};\\[3pt]
    m(n-k)+m+1-i, & \mbox{for $j=1$ and $m+1\leq i\leq m(n-k)+m$};\\[3pt]
    h_{(i,j-1)}', & \mbox{for $2\leq j\leq k$ and $1\leq i\leq m$}
    \end{array}\right.\label{eq-hooklengths-2}
\end{align}
where $h_{(i,j)}'$ denotes the hook-length of the cell $(i,j)$ in partition $((k-1)^m)$. Therefore,
\begin{align*}
f^{(k^m,1^{m(n-k)})}&=\frac{(mn)!}{\displaystyle{\prod_{1\leq i \leq m}}h_{(i,1)}\times \displaystyle{\prod_{m+1\leq i \leq m(n-k)+m}}h_{(i,1)}\times \displaystyle{\prod_{1\leq i \leq m, 2\leq j\leq k}}h_{(i,j)}}
\end{align*}
Substituting \eqref{eq-hooklengths-2} into the above formula, we obtain
\begin{align*}
f^{(k^m,1^{m(n-k)})}&=\frac{(mn)!}{[m(n-k)]!\times \displaystyle{\prod_{i=0}^{m-1}[m(n-k)+k+i]}\times
\displaystyle{\prod_{1\leq i \leq m,  2\leq j\leq k}}h'_{(i,j-1)}}.
\end{align*}
While, we see that
$$\prod_{1\leq i \leq m,  2\leq j\leq k}h'_{(i,j-1)}
=\prod_{(i,j)\in ((k-1)^m)}h'_{(i,j)}
=\frac{\displaystyle{\prod _{i=0}^{m-1}(k-1+i)!}}{\displaystyle{\prod_{0\leq i_1<i_2\leq m-1}(i_2-i_1)}},$$
where the second equality is obtained by applying \eqref{hook_length_schur_equa} to the partition $\lambda=((k-1)^m)$.
Thus, we have
\begin{align*}
f^{(k^m,1^{m(n-k)})}&=
\frac{(mn)!}{[m(n-k)]!\times\displaystyle{\prod_{i=0}^{m-1}[m(n-k)+k+i]}}
\times \frac{\displaystyle{\prod_{0\leq i_1<i_2\leq m-1}(i_2-i_1)}}{\displaystyle{\prod _{i=0}^{m-1}(k-1+i)!}}.
\end{align*}
Let
$$b_k=[m(n-k)]!\prod_{i=0}^{m-1}[m(n-k)+k+i]=\frac{(mn-mk)!(mn-mk+k-1)!}{(mn-mk+k+m-1)!}.$$
Then, for $1< k < n$, we have
\begin{align*}
\frac{(f^{(k^m,1^{m(n-k)})})^2}{f^{((k+1)^m,1^{m(n-k-1)})}f^{((k-1)^m,1^{m(n-k+1)})}}
=\frac{k+m-1}{k-1}\times
\frac{b_{k-1}b_{k+1}}{b_k^2}.
\end{align*}
Now it suffices to show that $b_k^2\leq b_{k-1}b_{k+1}$ for $1< k < n$.
Let $$b(z)=\frac{\Gamma(mn-mz+1)\Gamma(mn-mz+z)}{\Gamma(mn-mz+z+m)}$$ be the continuous function on $[1,n]$, where $\Gamma(z)$ is the Gamma function. Hence,  for $1\leq k \leq n$, we have $b_k=b(k)$, the value of $b(z)$ evaluated at $z=k$. To prove $b_k^2\leq b_{k-1}b_{k+1}$ for $1< k < n$, it suffices to show that $(\ln b(z))''\geq 0$
for $z\in [1,n]$. To this end, we first compute the logarithmic derivative of $b(z)$ as follows:
\begin{align*}
(\ln b(z))'
&=-m\psi(mn-mz+1)-(m-1)\psi(mn-mz+z)+(m-1)\psi(mn-mz+z+m)
\end{align*}
where $\psi(z)=(\ln \Gamma(z))'$ is the digamma function.
Then we obtain that
\begin{align*}
(\ln b(z))''=m^2\psi'(mn-mz+1)+(m-1)^2\psi'(mn-mz+z)-
(m-1)^2\psi'(mn-mz+z+m).
\end{align*}
It is known that $\psi'(z)=\sum_{k=0}^{\infty}\frac{1}{(z+k)^2}$ over $z\in (0,+\infty)$, and hence it is positive and decreasing, see \cite{abramowitz1972ny}. Thus, $(\ln b(z))''\geq 0$
for any $z\in [1,n]$. This completes the proof.
\qed

As an immediate consequence of Theorem \ref{thm-bona-lackner-sagan}, we obtain the following result, which shows that Theorem \ref{thm-bona} is true for any positive integer $m$.

\begin{cor}\label{general-bona} Suppose that $m,n$ are two positive integers.
\begin{itemize}
\item[(1)] For $\Lambda=\{(j^m,(n-j)^m) \mid \lceil\frac{n}{2}\rceil\leq j\leq n\}$, both $\{\ell^{\Lambda}_{mn,k}\}_{k=1}^{mn}$ and $\{i^{\Lambda}_{mn,k}\}_{k=1}^{mn}$ are log-concave.
\item[(2)] For $\Lambda=\{(j^m,1^{m(n-j)}) \mid 1\leq j\leq n\}$, both $\{\ell^{\Lambda}_{mn,k}\}_{k=1}^{mn}$ and $\{i^{\Lambda}_{mn,k}\}_{k=1}^{mn}$ are log-concave.
\end{itemize}
\end{cor}

\section{Proof by Schur positivity}\label{sec-schur_3}

The aim of this section is to give another proof of Theorem \ref{thm-bona-lackner-sagan} by using Schur positivity. Recall that a symmetric function is said to be Schur positive if it can be written as non-negative integer linear combination of Schur functions. By \eqref{exponential_specialization_s}, Theorem \ref{thm-bona-lackner-sagan} is implied by the following result.

\begin{thm}\label{special_shape1}
Suppose that $m$ and $n$ are two positive integers.
\begin{itemize}
\item[(1)] For $\lceil\frac{n}{2}\rceil< k < n$, the difference
\begin{align*}
s_{(k^m,(n-k)^m)}^2-s_{((k+1)^m,(n-k-1)^m)}s_{((k-1)^m,(n-k+1)^m)}
\end{align*}
is Schur positive.

\item[(2)] For $1 < k < n$, the difference
\begin{align*}
s_{(k^m,1^{m(n-k)})}^2-s_{((k+1)^m,1^{m(n-k-1)})}s_{((k-1)^m,1^{m(n-k+1)})}
\end{align*}
is Schur positive.
\end{itemize}
\end{thm}

Our proof of Theorem \ref{special_shape1} is based on some Schur positivity results due to Lam, Postnikov and Pylyavaskyy \cite{lam2007ajm}.
For vectors $v,w$ and a positive integer $n$, we assume that the operations $v+w$, $\frac{v}{n}$, $\lfloor v\rfloor$ and $\lceil v\rceil$ are performed coordinate-wise. In particular, we have well-defined operations $\lfloor \frac{\lambda+\mu}{2}\rfloor$ and $\lceil \frac{\lambda+\mu}{2}\rceil$ on pairs of any partitions.
If $\lambda,\mu$ are partitions with $\lambda_i\geq\mu_i$ for all $i\geq 1$, then the skew diagram $\lambda/\mu$ is the
diagram of $\lambda$ with the diagram of $\mu$ removed from its upper left-hand corner.
Lam, Postnikov and Pylyavaskyy obtained the following result,
which answered a conjecture of Okounkov \cite{okounkov1997aim}.

\begin{thm}[{\cite[Theorem 11]{lam2007ajm}}]\label{special_shape1_lemma}
 Given any two skew partitions $\lambda/\mu$ and $\nu/\rho$, the difference
$$s_{\lfloor \frac{\lambda+\nu}{2}\rfloor/
\lfloor \frac{\mu+\rho}{2}\rfloor}
s_{\lceil \frac{\lambda+\nu}{2}\rceil/
\lceil \frac{\mu+\rho}{2}\rceil}
-s_{\lambda/\mu}s_{\nu/\rho}$$
is Schur positive.
\end{thm}

 Given two partitions $\lambda$ and $\mu$, let $\lambda\cup\mu=(\nu_1,\nu_2,\nu_3,\ldots)$ be the partition obtained by rearranging all parts of $\lambda$ and $\mu$ in the weakly decreasing order. Let $\mathrm{sort}_1(\lambda,\mu):=(\nu_1,\nu_3,\nu_5,\ldots)$ and $\mathrm{sort}_2(\lambda,\mu):=(\nu_2,\nu_4,\nu_6,\ldots)$. Lam, Postnikov and Pylyavaskyy also obtained the following result, which was first conjectured by Fomin, Fulton, Li and Poon \cite{fomin2005ajm}.

\begin{thm}[{\cite[Corollary 12]{lam2007ajm}}]\label{special_shape2_lemma}
For any two partitions $\lambda$ and $\mu$, the
difference
$$s_{\mathrm{sort}_1(\lambda,\mu)}s_{\mathrm{sort}_2(\lambda,\mu)}-s_{\lambda}s_{\mu}$$ is Schur positive.
\end{thm}

We are now in the position to give a proof of Theorem \ref{special_shape1}.

\proof[Proof of Theorem \ref{special_shape1}.]
For $\lceil\frac{n}{2}\rceil< k < n$, taking  $\lambda=((k+1)^m,(n-k-1)^m)$, $\nu=((k-1)^m,(n-k+1)^m)$, and $\mu=\rho=\emptyset$ in Theorem \ref{special_shape1_lemma},
we obtain the Schur positivity of
$$s_{(k^m,(n-j)^m)}^2-s_{((k+1)^m,(n-k-1)^m)}s_{((k-1)^m,(n-k+1)^m)}.$$

For $1<k <n$, taking $\lambda=((k+1)^m,1^{m(n-k-1)})$, $\nu=((k-1)^m,1^{m(n-k+1)})$ and $\mu=\rho=\emptyset$ in Theorem \ref{special_shape1_lemma}, we obtain the Schur positivity of
\begin{align}\label{spp_1}
s_{(k^m,1^{m(n-k+1)})}s_{(k^m,1^{m(n-k-1)})}
-s_{((k+1)^m,1^{m(n-k-1)}}s_{((k-1)^m,1^{m(n-k+1)}}.
\end{align}
Taking $\lambda=(k^m,1^{m(n-k+1)})$ and $\mu=(k^m,1^{m(n-k-1)})$ in Theorem \ref{special_shape2_lemma}, we obtain the Schur positivity of
\begin{align}\label{spp_2}
s_{(k^m,1^{m(n-k)})}^2
-s_{(k^m,1^{m(n-k+1)})}s_{(k^m,1^{m(n-k-1)})}.
\end{align}
Combining \eqref{spp_1} and \eqref{spp_2}, we obtain the Schur positivity of
$$s_{(k^m,1^{m(n-k)})}^2-s_{((k+1)^m,1^{m(n-k-1)})}s_{((k-1)^m,1^{m(n-k+1)})}.$$
This completes the proof. \qed

As we mentioned at the end of Section \ref{sec-schur_2}, Theorem \ref{thm-bona} implies the log-concavity of certain sequences concerning longest increasing subsequences. The approach of this section to Theorem \ref{thm-bona} inspired us to study Conjecture \ref{chen_logcave_permu} and Conjecture \ref{bona_logcave_permu} from the viewpoint of Schur positivity. Note that for a fixed integer $n$, the Robinson-Schensted correspondence shows that
\begin{align}\label{RS_map_perm}
\ell_{n,k}{=}\sum_{\lambda\vdash n,\lambda_1=k}{\left(f^{\lambda}\right)}^2~~~
\mathrm{and}~~~
i_{n,k}{=}\sum_{\lambda\vdash n,\lambda_1=k}{f^{\lambda}}
\end{align}
for $1\leq k \leq n$. We have the following conjectures.

\begin{conj}\label{chen_logcave_permu_schur}
For
$1\leq k\leq n$, let $$ f_{n,k}=\sum_{\lambda\vdash n,\lambda_1=k}s_{\lambda}^2,$$
then $f_{n,k}^2-f_{n,k+1}f_{n,k-1}$ is Schur positive with the convention that $f_{n,0}=f_{n,n+1}=0$.
\end{conj}
\begin{conj}\label{bona_logcave_permu_schur}
For
$1\leq k\leq n$, let
$$g_{n,k}=\sum_{\lambda\vdash n,\lambda_1=k}s_{\lambda},$$ then $g_{n,k}^2-g_{n,k+1}g_{n,k-1}$ is Schur positive with the convention that $g_{n,0}=g_{n,n+1}=0$.
\end{conj}

It is readily to see that Conjecture \ref{chen_logcave_permu_schur} implies Conjecture \ref{chen_logcave_permu}, and Conjecture \ref{bona_logcave_permu_schur} implies Conjecture \ref{bona_logcave_permu} by \eqref{exponential_specialization_s}.
We have verified Conjecture  \ref{chen_logcave_permu_schur} for $n\leq 9$ and Conjecture \ref{bona_logcave_permu_schur} for $n \leq 20$.

Chen \cite{chen2008} also put forward some log-concavity conjecture about perfect matchings, which was turned into the form of Conjecture \ref{chen_logcave_perfectmatching} by  B\'{o}na, Lackner and Sagan.
For any fixed $n$, let
$\Theta$ be the set of partitions of $n$ all of whose column lengths are even.
Chen's conjecture can be stated as follows.

\begin{conj}[{\cite[Conjecture 1.5]{chen2008}}] \label{chen_logcave_perfectmatching}
For any fixed $n$, the sequence $\{i_{n,k}^{\Theta}\}_{k=1}^{n}$ is log-concave.
\end{conj}

Inspired by Conjectures \ref{chen_logcave_permu_schur} and \ref{bona_logcave_permu_schur}, we propose the following conjecture, which implies Conjecture \ref{chen_logcave_perfectmatching}.

\begin{conj}\label{chen_logcave_perfectmatching_schur}
For
$1\leq k\leq n$, let
$$g_{n,k}^{\Theta}=\sum_{\lambda\in \Theta,\lambda_1=k}s_{\lambda},$$ then $(g_{n,k}^{\Theta})^2-g_{n,k+1}^{\Theta}g_{n,k-1}^{\Theta}$ is Schur positive with the convention that $g_{n,0}^{\Theta}=g_{n,n+1}^{\Theta}=0$.
\end{conj}

This conjecture has been verified for $n\leq 30$. B\'{o}na, Lackner and Sagan further proposed a companion conjecture to Conjecture \ref{chen_logcave_perfectmatching}.

\begin{conj}[{\cite[Conjecture 4.3]{bona2016aoc}}]\label{bona_logcave_perfectmatching}
For any fixed $n$, the sequence $\{{\ell}_{n,k}^{\Theta}\}_{k=1}^{n}$ is log-concave.
\end{conj}

However, Conjecture \ref{bona_logcave_perfectmatching} does not admit a similar conjecture as Conjecture \ref{chen_logcave_perfectmatching_schur} as illustrated below. For $1\leq k\leq n$, let $$ f_{n,k}^{\Theta}=\sum_{\lambda \in \Theta,\lambda_1=k}s_{\lambda}^2.$$ In general, the difference
$(f_{n,k}^{\Theta})^2-f_{n,k+1}^{\Theta}f_{n,k-1}^{\Theta}$ is not Schur positive. For instance, when $n=10$, we have
\begin{align*}
&f^{\Theta}_{10,2}=s_{(2,2,2,2,1,1)}^2+s_{(2,2,1,1,1,1,1,1)}^2,\\
&f^{\Theta}_{10,3}=s_{(3,3,2,2)}^2+s_{(3,3,1,1,1,1)}^2,\\
&f^{\Theta}_{10,4}=s_{(4,4,1,1)}^2.
\end{align*}
However, the symmetric function
$(f^{\Theta}_{10,3})^2-f^{\Theta}_{10,2}f^{\Theta}_{10,4}$ is not Schur positive by computer exploration using the open-source
mathematical software \texttt{Sage}~\cite{sage} and its algebraic
combinatorics features developed by the \texttt{Sage-Combinat}
community~\cite{Sage-Combinat}.

\noindent{\bf Acknowledgements.} This work was supported by the National Science Foundation of China.

\end{document}